\documentclass[12pt]{amsart}
\usepackage{bm}
 \usepackage{graphicx}

 \newtheorem{thm}{Theorem}[section]
 \newtheorem{cor}[thm]{Corollary}
 
 \newtheorem{prop}[thm]{Proposition}
 \theoremstyle{definition}
 
 \theoremstyle{remark}
 
 \numberwithin{equation}{section}

 \newcommand{\Real}{\mathbb{R}}

\begin{document}

\title[On surfaces in three dimensional contact manifolds]
{On surfaces in three dimensional contact manifolds}

\author{Paul W.Y. Lee}
\email{wylee@math.cuhk.edu.hk}
\address{Room 216, Lady Shaw Building, The Chinese University of Hong Kong, Shatin, Hong Kong}

\date{\today}

\begin{abstract}
In this paper, we introduce two notions on a surface in a contact manifold. The first one is called degree of transversality (DOT) which measures the transversality between the tangent spaces of a surface and the contact planes. The second quantity, called curvature of transversality (COT), is designed to give a comparison principle for DOT along characteristic curves under bounds on COT. In particular, this gives estimates on lengths of characteristic curves assuming COT is bounded below by a positive constant.

We show that surfaces with constant COT exist and we classify all graphs in the Heisenberg group with vanishing COT. This is accomplished by showing that the equation for graphs with zero COT can be decomposed into two first order PDEs, one of which is the backward invisicid Burgers' equation. Finally we show that the p-minimal graph equation in the Heisenberg group also has such a decomposition. Moreover, we can use this decomposition to write down an explicit formula of a solution near a regular point.
\end{abstract}

\maketitle

\section{Introduction}

Motivated by the isoperimetric problems in three-dimensional contact manifolds or pseudo-hermitian geometry, surfaces in these manifolds have received an increasing interest in recent years (e.g. \cite{ChHwMaYa1, ChHwMaYa2,DaGaNhPa,HuRiRo,RiRo}). These surfaces are foliated by curves, called characteristic curves, which play a very important role in the understanding of submanifold geometry in these spaces. Recall that the tangent spaces of a surface intersect transversely with the three dimensional contact planes at generic points. This defines a line field on the surface and the leaves of the corresponding foliation are the characteristic curves. In this paper, we study these curves from the point of view of comparison geometry. We introduce two quantities which are closely related to these characteristic curves. The first quantity is called degree of transversality (DOT) which measures how transverse the intersections are between the tangent spaces of the surface and the contact planes. In particular, DOT is infinite at a point if the tangent plane coincides with the contact plane there. Such a point is called singular.

It is pointed out by the referee that DOT should be closely related to the imaginary curvature introduced by \cite{ArFe1,ArFe2}. In fact, we show that DOT and the imaginary curvature coincide up to a sign in Proposition \ref{imaginarycurv}. Imaginary curvature is a quantity defined on a surface in the simplest subriemannian manifold, the Heisenberg group. It was used in \cite{ArFe1} to study subriemannian analogue of focal points of a surface and in \cite{ArFe2} to study the horizontal Hessian of the subriemannian distance in the Heisenberg group. In \cite{AgLe2} motivated by the earlier work \cite{AgLe1}, a version of subriemannian Hessian was introduced. One can generalize the above mentioned results to all three dimensional Sasakian manifolds using DOT and this will be reported in a forthcoming work.

Next, we introduce the second quantity, called curvature of transversality (COT), which gives a comparison principle for DOT along characteristic curves. This comparison principle is similar to the one for Ricci curvature which compares Jacobi fields along geodesics. In particular, we show that characteristic curves must hit two singular points if COT is bounded below by a positive constant and they hit at most one singular point if COT is bounded above by a negative constant. All these are accomplished in Section \ref{submanifold}.

The next natural question is whether there are surfaces with constant COT. In Section \ref{COTeg}, we show that the answer to this question is positive. The main result of this paper is the classification of graphs in the Heisenberg group with vanishing COT. This is done in Section \ref{zeroCOT}.

The second order PDE (\ref{zCOT}) satisfied by graphs  with vanishing COT in the Heisenberg group is very similar to the p-minimal graph equation studied in \cite{Pa1,Pa2,GaPa,ChHwMaYa1}. Recall that a surface is minimal if the mean curvature vanishes everywhere. The Berstein theorem says that the graph of a function in $\Real^2$ is a minimal surface in $\Real^3$ if and only if the function is linear. In the subriemannian case, there is an analogue of the mean curvature called $p$-mean curvature and a surface is $p$-minimal if $p$-mean curvature vanishes. The graph of a function $f$ in $\Real^2$ is a $p$-minimal surface in the Heisenberg group if it satisfies the $p$-minimal graph equation (\ref{pminimal}). The following families of global solutions to the $p$-minimal graph equation were found in \cite{Pa2}:
\begin{equation}\label{Bern}
\begin{split}
&f(x,y)=ax+by+c,\\
& f(x,y)=-abx^2 + (a^2-b^2)xy + aby^2 + g(-bx + ay)
\end{split}
\end{equation}
where $a,b,c$ are constants and $g:\Real\to\Real$ is any $C^2$ function.

The subriemannian version of Berstein theorem proved by \cite{GaPa,ChHwMaYa1} says that any $C^2$ solution of the $p$-minimal graph equation is given by (\ref{Bern}).

In this paper, we show that the PDE (\ref{zCOT}) satisfied by graphs with vanishing COT in the Heisenberg group and the $p$-minimal graph equation (\ref{pminimal}) can be split into two first order PDE, one of them is the inviscid Burgers' equation. By using the method of characteristics, we obtain an explicit formula for any local solution to the equations near a regular point. Using this local solution of the $p$-minimal graph equation, we recover the global solutions (\ref{Bern}) in a simple way. All these will be accomplished in the last two sections of the paper.

\smallskip

\section*{Acknowledgment}
The author would like to thank the referee for several very constructive comments which improve the exposition of the paper and stimulate further research on the subject.

\smallskip

\section{Submanifolds in Contact Geometry}\label{submanifold}

In this section, we recall and introduce several notions on submanifolds in contact geometry which are needed in this paper. Let us start with some basic notions in contact geometry. Recall that a three dimensional manifold $M$ is \textit{contact} if there is a 1-form $\alpha_0$, called a \textit{contact form}, such that $d\alpha_0$ is non-degenerate (i.e. the map $v\mapsto d\alpha_0(v,\cdot)$ is a bijection). The kernel of the contact form $\alpha_0$ defines a distribution $\Delta$ (a vector subbundle of the tangent bundle). Note that if $\alpha_0$ is a contact form, then so is $f\alpha_0$ where $f$ is a nonzero function on $M$. If we fix a smoothly varying inner product $\left<\cdot,\cdot\right>$, called a \textit{subriemannian metric}, on the distribution $\Delta$, then there is a unique contact form $\alpha_0$ such that the restriction of the non-degenerate form $d\alpha_0$ to the distribution $\Delta$ coincides with the volume form induced by the subriemannian metric (i.e. $d\alpha_0(v_1,v_2)=1$, where $\{v_1,v_2\}$ is orthonormal with respect to the subriemannian metric). We also define a vector field $v_0$, called the \textit{Reeb field}, by the conditions $\alpha_0(v_0)=1$ and $d\alpha_0(v_0,\cdot)=0$.

Let $N$ be a submanifold. A point $x$ on $N$ is \textit{regular} if $T_xN$ and $\Delta_x$ intersect transversely. Otherwise it is called \textit{singular}. For each regular point $x$ on the submanifold $N$, we can pick an \textit{adapted frame} $\{v_1,v_2\}$ around $x$. This is a pair of vector fields $v_1$ and $v_2$ of the ambient manifold $M$ which satisfy three conditions. First, they form an orthonormal basis of the distribution $\Delta$ with respect to the given subriemannian metric. Second, the first vector field $v_1$ is contained in the line field $TN\cap\Delta$. Third, they are oriented by the condition $d\alpha_0(v_1,v_2)=1$.  Note that if $v_1,v_2$ is such a frame, then $-v_1,-v_2$ is the only other choice. The integral curves of $v_1$ are tangent to the line field $TN\cap\Delta$ and they are called \textit{characteristic curves}.

Next, we introduce the degree of transversality (DOT) mentioned in the introduction. For this, let us fix a characteristic curve $\gamma$. For each point $x$ on the characteristic curve, we fix an adapted frame $v_1,v_2$. Since $x$ is a regular point, $v_2$ is not contained in the tangent bundle $TN$ of the submanifold $N$. Therefore, there is a function $\mathfrak a$ such that $v_0-\mathfrak av_2$ is contained in $TN$. We call the function $\mathfrak a$ defined along the characteristic curve by \textit{degree of transversality} (DOT). Note that $\mathfrak a$ approaches $\pm\infty$ as the points on the characteristic curve approaches a singular point. Note also that DOT depends only on the orientation of the characteristic curve. If we pick the opposite orientation, then $-v_1,-v_2$ is another adapted frame and DOT is given by $-\mathfrak a$ in this case. In particular, $|\mathfrak a|$ is well-defined on the set of all regular points.

Next, we give a simple expression for $|\mathfrak a|$.

\begin{prop}\label{imaginarycurv}
Assume that $N$ is given by a level of a smooth function $g$ and let $\nabla_H g$ be the horizontal gradient defined by
\[
\nabla _H g=(v_1g) v_1+(v_2g) v_2.
\]
Then $|\mathfrak a|$ satisfies
\[
|\mathfrak a|=\frac{|v_0g|}{|\nabla_Hg|}.
\]
\end{prop}

\begin{proof}
By assumptions, both $v_1$ and $v_0-\mathfrak a v_2$ are contained in the tangent bundle $TN$. Therefore, $v_1g=v_0g-\mathfrak a v_2g=0$. Hence, $|\mathfrak a|=\frac{|v_0g|}{|v_2g|}=\frac{|v_0g|}{|\nabla_Hg|}$.
\end{proof}

In the Heisenberg group, the quantity $\frac{v_0g}{|\nabla_Hg|}$ is the imaginary curvature introduced in \cite{ArFe1,ArFe2}. Proposition \ref{imaginarycurv} shows a close relation between this and DOT.

The \textit{curvature of transversality} (COT) is defined by
\[
\mathfrak r=v_1\mathfrak a-\mathfrak a^2.
\]
Note that COT is defined on the set of all regular points of $N$, not just along a characteristic curve.

By design, if $\gamma(\cdot)$ is a characteristic curve, then the following equation holds
\begin{equation}\label{Riccati}
\frac{d}{dt}\mathfrak a(\gamma(t))=\mathfrak a(\gamma(t))^2+\mathfrak r(\gamma(t)).
\end{equation}

Therefore, DOT and COT satisfy the following comparison principle.

\begin{prop}\label{compare}
Let $\gamma$ be a characteristic curve of $v_1$. Assume that $\mathfrak r(\gamma(t))\leq k(t)$ (resp. $\geq$, $<$, $>$) and let $c(t)$ be the solution of
\[
\dot c(t)=c(t)^2+k_2(t)
\]
with $c(0)=\mathfrak a(\gamma(0))$. Then
\[
\mathfrak a(\gamma(t))\leq c(t)\quad (\text{resp. } \geq, <, >)
\]
for all $t\geq 0$.
\[
\mathfrak a(\gamma(t))\geq c(t)\quad (\text{resp. } \leq, >, <)
\]
for all $t\leq 0$.
\end{prop}

\begin{proof}
The difference $d(t):=\mathfrak a(\gamma(t))-c(t)$ satisfies
\[
\begin{split}
\dot d&=\dot{\mathfrak a}-\dot c\\
&\leq\mathfrak a(\gamma)^2-c^2\quad (\text{resp.} \geq, <, >)\\
&=(\mathfrak a(\gamma)+c)\,d.
\end{split}
\]

Therefore, by Gronwall's inequality, $\mathfrak a(\gamma(t))\leq c(t)$ (resp. $\geq$, $<$, $>$) for all $t\geq 0$. Similar arguments show that $\mathfrak a(\gamma(t))\geq c(t)$ (resp. $\leq$, $>$, $<$) if $t<0$.
\end{proof}

When COT is constant, DOT can be computed explicitly (see section \ref{COTeg} for examples of surfaces with constant COT). If we combine this with Proposition \ref{compare}, then we obtain the following.

\begin{prop}
Let $\gamma$ be a characteristic curve tangent to the vector field $v_1$. Assume that $\mathfrak r(\gamma(t))\leq k$ (resp. $\geq$, $<$, $>$) for some constant $k$. Then the following holds
\[
\begin{split}
&\mathfrak a(t) \leq (\text{resp.} \geq, <, >) \\
&\begin{cases}
\frac{\sqrt k(\cos(t\sqrt k)\mathfrak a(\gamma(0))+\sqrt k\sin(t\sqrt k))}{-\sin(t\sqrt k)\mathfrak a(\gamma(0))+\sqrt k\cos(t\sqrt k)}
 & \mbox{if } k>0\\
\frac{\mathfrak a(\gamma(0))}{1-\mathfrak a(\gamma(0))t} & \mbox{if } k=0\\
\frac{\sqrt{-k}(\cosh(t\sqrt{-k})\mathfrak a(\gamma(0))-\sqrt{-k}\sinh(t\sqrt{-k}))}{-\sinh(t\sqrt{-k})\mathfrak a(\gamma(0)) +\sqrt{ -k}\cosh(t\sqrt{-k})} & \mbox{if } k<0.
\end{cases}
\end{split}
\]
\end{prop}

As a consequence, we have the following results on the singular set.

\begin{cor}
Let $\gamma$ be a characteristic curve of $v_1$.
\begin{itemize}
\item If $\mathfrak r(\gamma(t))\leq k\leq 0$ for all time $t$ and $|\mathfrak a(\gamma(0))|\leq \sqrt{-k}$, then there is no singular point along $\gamma$.
\item If $\mathfrak r(\gamma(t))\leq k\leq 0$  for all time $t$ and $|\mathfrak a(\gamma(0))|> \sqrt{-k}$, then there is at most one singular point along $\gamma$.
\end{itemize}

Let $t_{-\infty}<0$ and $t_{\infty}>0$ be the time such that $\gamma(t_{-\infty})$ and $\gamma(t_{\infty})$ are singular points.

\begin{itemize}
\item If $\mathfrak r(\gamma(t))\geq k\geq 0$ for all time $t\geq 0$, then
\[
t_\infty\leq
\begin{cases}
\frac{1}{\sqrt k}\cot^{-1}\left(\frac{\mathfrak a(\gamma(0))}{\sqrt k}\right) & \mbox{if } k>0\\
\frac{1}{\sqrt{-k}}\coth^{-1}(\mathfrak a(\gamma(0))\sqrt{-k}) & \mbox{if } k<0,\mathfrak a(\gamma(0))>\sqrt{-k}\\
\frac{1}{\mathfrak a(\gamma(0))} & \mbox{if } k=0, \mathfrak a(\gamma(0))> 0.
\end{cases}
\]
Moreover, equality holds only if  $\mathfrak r(\gamma(t))=k$ for all $t\geq 0$.

\item If $\mathfrak r(\gamma(t))\geq k$ for all time $t\leq 0$, then
\[
t_{-\infty}\leq
\begin{cases}
-\frac{\pi}{\sqrt k}+\frac{1}{\sqrt k}\cot^{-1}\left(\frac{\mathfrak a(\gamma(0))}{\sqrt k}\right) & \mbox{if } k>0\\
\frac{1}{\sqrt{-k}}\coth^{-1}(\mathfrak a(\gamma(0))\sqrt{-k}) & \mbox{if } k<0,\mathfrak a(\gamma(0))<-\sqrt{-k}\\
\frac{1}{\mathfrak a(\gamma(0))} & \mbox{if } k=0, \mathfrak a(\gamma(0))< 0.
\end{cases}
\]
Moreover, equality holds only if  $\mathfrak r(\gamma(t))=k$ for all $t\leq 0$.

\item If $\mathfrak r(\gamma(t))\geq k>0$ for all time $t$, then there are two singular points along $\gamma$ and the length of the characteristic curve $\gamma$ is at most $\frac{\pi}{\sqrt k}$. Moreover, equality holds only if $\mathfrak r(\gamma(t))=k$ for all $t$.

\item In particular, if $\mathfrak r\geq k> 0$ on the submanifold $N$ and the singular set is bounded, then $N$ is compact.
\end{itemize}
\end{cor}

\smallskip

\section{Examples of Submanifolds with Constant COT}\label{COTeg}

In this section, we show that surfaces with constant COT exist. Let us first give another characterization of COT. Let $N$ be a given submanifold of a three dimensional contact manifold. Let $v_1,v_2$ be an adapted frame and $v_0$ be the Reeb field. Let $a_{ij}^k:N\to\Real$ be the structure constants defined by

\begin{equation}\label{bracket}
[v_i,v_j]=\sum_{k=0}^2a_{ij}^kv_k,
\end{equation}
where $i,j=0,1,2$.

\begin{prop}\label{COTformula}
Under the notations introduced above, COT is given by
\[
\mathfrak r=-a_{01}^2-\mathfrak aa_{12}^2.
\]
\end{prop}

\begin{proof}
We consider the dual version of (\ref{bracket}),
\[
d\alpha_k=-\sum_{0\leq i<j\leq 2}a_{ij}^k\alpha_i\wedge\alpha_j.
\]

By the definition of adapted frame, we have $-a_{12}^0=d\alpha_0(v_1,v_2)=1$. Therefore,
\[
d\alpha_0=\alpha_1\wedge\alpha_2-a_{01}^0\alpha_0\wedge\alpha_1-a_{02}^0\alpha_0\wedge\alpha_2.
\]

By the definition of the Reeb field $v_0$, we also have $d\alpha_0(v_0,\cdot)=0$ and so
$a_{01}^0=a_{02}^0=0$.

The two vector fields $v_1$ and $v_0-\mathfrak a v_2$ are tangent to the submanifold $N$. Therefore, the bracket
\[
\begin{split}
&[v_1,v_0-\mathfrak a v_2]\\
&=-a_{01}^1v_1-a_{01}^2v_2-(v_1\mathfrak a)v_2-\mathfrak a(-v_0+a_{12}^1v_1+a_{12}^2v_2)\\
&=-(a_{01}^1+\mathfrak aa_{12}^1)v_1+\mathfrak a (v_0-\mathfrak a v_2)+(\mathfrak a^2-a_{01}^2-(v_1\mathfrak a)-\mathfrak aa_{12}^2)v_2.
\end{split}
\]
is also tangent to $N$. Hence, $[v_1,v_0-\mathfrak a v_2]$ is a linear combination of $v_1$ and $v_0-\mathfrak a v_2$ and the following holds on $N$
\[
\mathfrak r=v_1\mathfrak a-\mathfrak a^2=-a_{01}^2-\mathfrak aa_{12}^2.
\]
\end{proof}

Using Proposition \ref{COTformula}, it is not hard to construct examples with constant COT. Recall that $SU(2)$, the special unitary group, consists of $2\times2$ matrices with complex coefficients and determinant 1. The Lie algebra $su(2)$ consists of skew Hermitian matrices with trace zero. The following two elements in $su(2)$
\[
v_1=\left(
\begin{array}{cc}
0 & 1/2\\
-1/2 & 0
\end{array}\right),\quad v_2=\left(
\begin{array}{cc}
0 & i/2\\
i/2 & 0
\end{array}\right)
\]
defines a subriemannian structure on $SU(2)$. The distribution is given by the span of the two left-invariant vector fields corresponding to $v_1$ and $v_2$. The subriemannian metric is the one which satisfies  $\left<v_i,v_j\right>=\delta_{ij}$, $i=1,2$. The Reeb field $v_0$, in this case,  is given by
\[
v_0=\left(
\begin{array}{cc}
-i/2 & 0\\
0 & i/2
\end{array}\right).
\]

Let $N$ be a surface which is foliated by the integral curves of the left invariant vector field defined by $v_1$. Then $v_1,v_2$ is an adapted frame and we have $a_{12}^2=0$ and $a_{01}^2=-1$. Therefore, by Proposition \ref{COTformula}, $\mathfrak r\equiv 1$ on $N$.

A specific example is given by the image of the following map.
\[
(\theta_1,\theta_2)\mapsto\left(\begin{array}{cc}
\cos(\theta_1/2) & \sin(\theta_1/2)\\
-\sin(\theta_1/2) & \cos(\theta_1/2)\\
\end{array}\right)
\left(\begin{array}{cc}
\cos(\theta_2/2) & i\sin(\theta_2/2)\\
i\sin(\theta_2/2) & \cos(\theta_2/2)\\
\end{array}\right).
\]
By rescaling the subriemannian structure, we can obtain examples with any positive constant COT.

For surfaces with constant negative COT, we consider the special linear group $SL(2)$, the set of all $2\times 2$ matrices with real coefficients and determinant 1. The Lie algebra $sl(2)$ is the set of all $2\times 2$ real matrices with trace zero. The left invariant vector fields of the following two elements in $sl(2)$
\[
v_1=\left(
\begin{array}{cc}
1/2 & 0\\
0 & -1/2
\end{array}\right),\quad v_2=\left(
\begin{array}{cc}
0 & 1/2\\
1/2 & 0
\end{array}\right)
\]
span a distribution $\Delta$ on $SL(2)$. The subriemannian metric on $SL(2)$ is defined by $\left<v_i,v_j\right>=\delta_{ij}$, $i=1,2$. The Reeb field in this case is $v_0$, where
\[
v_0=\left(
\begin{array}{cc}
0 & -1/2\\
1/2 & 0
\end{array}\right).
\]

Similar to the case in $SU(2)$, surfaces $N$ foliated by integral curves of $v_1$ will have constant COT. Since  $a_{12}^2=0$ and $a_{01}^2=1$, $N$ will have COT equal to $-1$.

For surfaces with constant zero COT, see section \ref{zeroCOT}.

\smallskip

\section{Graphs with Vanishing COT in the Heisenberg Group}\label{zeroCOT}

In this section, we consider graphs over the $xy$-plane in the Heisenberg group. Recall that the Heisenberg group is a subriemannian manifold on $\Real^3$ with distribution $\Delta$ spanned by two vector fields $u_1=\partial_x-\frac{1}{2}y\partial_z$ and $u_2=\partial_y+\frac{1}{2}x\partial_z$. The subriemannian metric is defined by declaring that $u_1$ and $u_2$ are orthonormal.

Let $f:U\to\Real$ be a function defined on a domain $U$ and let $N$ be the graph of $f$
\begin{equation}\label{graph}
N=\{(x,y,f(x,y))\in\Real^3|x,y\in\Real\}.
\end{equation}

In this case, the Reeb field $v_0$ is $-\partial_z$ and we can choose the adapted frame by
\[
v_1=\frac{(x-2f_y)\partial_x+(y+2f_x)\partial_y+(xf_x+yf_y)\partial_z}{\sqrt{(x-2f_y)^2+(y+2f_x)^2}}
\]
and
\[
v_2=\frac{(-y-2f_x)\partial_x+(x-2f_y)\partial_y+\frac{1}{2}(y^2+2yf_x+x^2-2xf_y)\partial_z}{\sqrt{(x-2f_y)^2+(y+2f_x)^2}}.
\]

A computation using the definition of DOT and COT shows that
\begin{equation}\label{GDOT}
\mathfrak a=-\frac{2}{\sqrt{(x-2f_y)^2+(y+2f_x)^2}}
\end{equation}
and
\begin{equation}\label{GCOT}
\begin{split}
\mathfrak r =&\frac{4(x-2f_y)(y+2f_x)(f_{yy}-f_{xx})+2(1-2f_{xy})(y+2f_x)^2}{((x-2f_y)^2+(y+2f_x)^2)^2}\\
&+\frac{2(x-2f_y)^2 (1+2f_{xy})}{((x-2f_y)^2+(y+2f_x)^2)^2}.
\end{split}
\end{equation}

In this section, we consider the following equation satisfied by graphs with zero COT in the Heisenberg group
\begin{equation}\label{zCOT}
\begin{split}
&2(x-2f_y)(y+2f_x)(f_{yy}-f_{xx})\\
&+(1-2f_{xy})(y+2f_x)^2 +(x-2f_y)^2 (1+2f_{xy})=0.
\end{split}
\end{equation}

A point $(x_0,y_0)$ is regular if and only if DOT is finite. So either $x-2f_y\neq 0$ or $y+2f_x\neq 0$. A computation shows the following.

\begin{prop}\label{split}
Let $U$ be an open set where $x-2f_y\neq 0$ (resp. $y+2f_x\neq 0$). Let $f$ be a $C^2$ solution of the equation (\ref{zCOT}) on $U$. Then
\[
g(x,y)=\frac{y+2f_x}{x-2f_y}\quad \left(\text{resp. } h(x,y)=\frac{x-2f_y}{y+2f_x}\right)
\]
is a solution of the backward inviscid Burgers' equation
\begin{equation}\label{Burgers}
g_y=gg_x \quad \left(\text{resp. } h_x=hh_y\right) .
\end{equation}
\end{prop}

Proposition \ref{split} shows that the equation (\ref{zCOT}) splits (near a regular point) into two first order PDEs, one of which is the backward inviscid Burgers' equation. Therefore, all solutions of (\ref{zCOT}) can be found near a regular point by first solving (\ref{Burgers}) by the method of characteristics (see for instance \cite{Ev}). Then substituting the resulting solution $g$ (resp. $h$) into
\[
y+2f_x=(x-2f_y)g\quad (\text{resp. } x-2f_y=(y+2f_x)h)
\]
and applying the method of characteristics again to obtain a solution $f$ of (\ref{zCOT}). Next we apply this observation to the proof of the main result. First, we have the following result on the singular set of any solution of (\ref{zCOT}). From now on, we denote a singular point $(x_0,y_0,f(x_0,y_0))$ on the graph $N$ simply by $(x_0,y_0)$.

\begin{prop}\label{C1curve}
Let $f$ be a (local) $C^2$ solution of (\ref{zCOT}) and let $N$ be the graph defined by (\ref{graph}). Then, for each singular point $(x_0,y_0)$ on $N$, there is a neighborhood around $(x_0,y_0)$ on which the singular set is given by a $C^1$ curve. Moreover, this $C^1$ curve is defined either by the equation $x-2f_y=0$ or $y+2f_x=0$. In particular, there is no isolated singular point on $N$.
\end{prop}

\begin{proof}
First, note that $\partial_x(x-2f_y)=1-2f_{xy}$ and $\partial_y(y+2f_x)=1+2f_{xy}$. Therefore, either $x-2f_y=0$ or $y+2f_x=0$ is a $C^1$ curve in a neighborhood of $(x_0,y_0)$. Without loss of generality, assume that $1-2f_{xy}\neq 0$ in a neighborhood of $(x_0,y_0)$ and we let $\Gamma$ be the $C^1$ curve defined by $x-2f_y=0$. The singular set is contained in $\Gamma$ and a point in $\Gamma$ is a singular if and only if $y+2f_x=0$ at that point.

Let $(x_i,y_i)$ be a sequence of regular points on $\Gamma$ which converges to $(x_0,y_0)$ as $i$ goes to $\infty$. By (\ref{zCOT}), we have
\[
(1-2f_{xy}(x_i,y_i))(y_i+2f_x(x_i,y_i))^2=0.
\]

Since $y_i+2f_x(x_i,y_i)\neq 0$, we must have $1-2f_{xy}(x_i,y_i)=0$ for all $i$. If we let $i$ goes to $\infty$, then we obtain $1-2f_{xy}(x_0,y_0)=0$ which is a contradiction. Therefore, there must be a neighborhood of the point $(x_0,y_0)$ in $\Gamma$ consists only of singular points.
\end{proof}

Next, we show that the domain of the function $f$ is foliated by lines where the functions $g$ and $h$ are constant.

\begin{thm}\label{constantg}
Let $f$ be a $C^2$ (local) solution to the equation (\ref{zCOT}) and let $g$ (resp. $h$) be as in Proposition \ref{split}. Then the domain of $f$ is foliated by lines and the function $g$ (resp. $h$) is constant or infinite along these lines. If a point $(a,b)$ satisfies $a-f_y(a,b)\neq 0$ (resp. $b+f_x(a,b)\neq 0$), then the line which passes through the point $(a,b)$ is given by
\[
x=-g(a,b)(y-b)+a \quad (\text{resp. } y=-h(a,b)(x-a)+b).
\]
\end{thm}

\begin{proof}
We only proof the statement for $g$ only. The one for $h$, being very similar, will be omitted. On the set where $x-f_y(x,y)\neq 0$, we can define a curve $\gamma(\cdot)$ by $\dot\gamma(t)=(-g(\gamma(t)),1)$ and $\gamma(0)=(a,b)$. Since
\[
\frac{d}{dt}g(\gamma(t))=-g_xg+g_y=0,
\]
it follows that the curve $\gamma$ is a straight line given by
\[
\gamma(t)=(-g(a,b)\,t+a,t+b)
\]
and $g$ is constant along $\gamma$.

On the set of regular points where $x-f_y(x,y)=0$, $g$ is infinite and $h$ is zero. The same argument as above shows that this set is foliated by horizontal lines. On these lines, $h$ vanishes and, therefore, $g$ is infinite. Therefore, it remains to consider what happen around a singular point.

Let $(x_0,y_0)$ be a singular point. We follow the notations in the proof of Proposition \ref{C1curve} and assume, without loss of generality, that $1-2f_{xy}\neq 0$ in a neighborhood $U$ of $(x_0,y_0)$. We let $\Gamma$ be the $C^1$ curve defined by $x-2f_y=0$. Let $(x_i,y_i)$ be a sequence of point outside $\Gamma$ which converges to $(x_0,y_0)$. By (\ref{zCOT}), we must have
\[
g^2+\frac{2(f_{yy}-f_{xx})g}{1-2f_{xy}}+\frac{1+2f_{xy}}{1-2f_{xy}}=0
\]
along $(x_i,y_i)$.

If we let $i$ goes to $\infty$, then we see that $g(x_i,y_i)$ can converge to at most two finite values. On the other hand, we have $x_0-2f_y(x_0,y_0)=0=y_0+2f_x(x_0,y_0)$. It follows that
\begin{equation}\label{ftaylor}
\begin{split}
f(x,y)=&f(x_0,y_0)-\frac{y_0}{2}(x-x_0)+\frac{x_0}{2}(y-y_0)\\
&+a_1(x-x_0)^2+a_2(x-x_0)(y-y_0)\\
&+a_3(y-y_0)^2+o(|x-x_0|^2+|y-y_0|^2).
\end{split}
\end{equation}

If we substitute (\ref{ftaylor}) into the definition of $g$, we obtain
\[
g(x,y)=\frac{4a_1(x-x_0)+(1+2a_2)(y-y_0)+o(\sqrt{(x-x_0)^2+(y-y_0)^2})}{(1-2a_2)(x-x_0)-4a_3(y-y_0)+o(\sqrt{(x-x_0)^2+(y-y_0)^2})}.
\]

Let $m_2(y-y_0)=m_1(x-x_0)$. Then the above equation of $g$ becomes
\[
\frac{4a_1m_2+(1+2a_2)m_1+o(1)}{(1-2a_2)m_2-4a_3m_1+o(1)}.
\]
It follows that $g$ approaches the same value along any line (since otherwise it will approach infinitely many values). This gives an extension of $g$ to the whole neighborhood $U$ and we denote this extension again by $g$.

By the earlier discussion, the vector field $(x,y)\mapsto (-g(x,y),1)$ is $C^1$ outside the set $\Gamma$ and the integral curves are straight lines. Let us denote the line
\begin{equation}\label{line}
x=-g(a,b)(y-b)+a
\end{equation}
corresponding to the point $(a,b)$ by $l_{(a,b)}$ (Note that $l_{(a,b)}$ still make sense even if $(a,b)$ is singular). First, let us assume that $l_{(x_0,y_0)}$ intersect $\Gamma$ transversely and show that $g$ is constant on $l_{(x_0,y_0)}$.

By the discussion above, it is enough to show that $g$ is constant in a neighborhood of $(x_0,y_0)$ inside the line $l_{(x_0,y_0)}$. Suppose it is not true. Then  we can find a sequence of points $(x_i,y_i)$ on $l_{(x_0,y_0)}$ converging to $(x_0,y_0)$ such that $g(x_i,y_i)\neq g(x_0,y_0)$ for all $i$. Let us fix $j$ and consider the line $l_{(x_j,y_j)}$. Since $l_{(x_0,y_0)}$ and $l_{(x_j,y_j)}$ intersect transversely, $l_{(x_i,y_i)}$ and $l_{(x_j,y_j)}$ intersect transversely as well for $i$ large enough. Moreover, we can assume that the intersections are close to $(x_j,y_j)$ and so they are contained in $U$. By choosing $j$ large enough, we can also assume that the line segment between these intersections and the points $(x_i,y_i)$ consist only of regular points. But then, by the method of characteristics discussed above $g(x_j,y_j)=g(x_i,y_i)$ for all $i$ large enough which is a contradiction. It follows that $g$ is constant along $l_{(x_0,y_0)}$ if it intersects $\Gamma$ transversely.

Next, we assume that $l_{(x_0,y_0)}$ intersects $\Gamma$ tangentially and there is a sequence of points $(x_i,y_i)$ on $l_{(x_0,y_0)}$ such that the lines $l_{(x_i,y_i)}$ either intersect $\Gamma$ transversely or do not intersect it at all for all $i$. The previous claim shows that $g$ is constant along each line $l_{(x_i,y_i)}$. Therefore, by using the previous argument, we also see that $g$ is constant along $l_{(x_0,y_0)}$ as well in this case.

Finally, it remains to consider the case where $l_{(x,y)}$ is tangent to $\Gamma$ for all points $(x,y)$ in a neighborhood of $(x_0,y_0)$ inside $\Gamma$. By the proof of Proposition \ref{C1curve}, $\Gamma$ can be parametrized by the $C^1$ path $t\mapsto (\varphi(t),t)$. It follows that $\dot\varphi(t)=-g(\varphi(t),t)$. A computation shows that $\frac{d}{dt}g(\varphi(t),t)=0$. Therefore, $\Gamma$ is everywhere tangent to lines of the same slope. Therefore, $\Gamma$ is a straight line and $g$ is constant along $\Gamma$.
\end{proof}

We call a $C^2$ solution $f$ of the equation (\ref{zCOT}) entire if $f$ is defined everywhere on the xy-plane. As a consequence of the above theorem, the functions $g$ and $h$ are constant functions if $f$ is entire.

\begin{cor}\label{constantgcor}
Let $f$ be an entire solution of the equation (\ref{zCOT}). Then the functions $g$ and $h$ defined in Proposition \ref{split} are constants.
\end{cor}

\begin{proof}
If $g$ is different at two points, then the lines corresponding to these two points given by (\ref{line}) have different slopes. Hence they must intersect. But this contradicts Theorem \ref{constantg}.
\end{proof}

Finally we prove the classification result mentioned in the introduction.

\begin{thm}\label{main}
Let $f$ be an entire solution of the equation (\ref{zCOT}). Then there are constants $c_1$, $c_2$, and a function $F$ such that the solution $f$ is given locally in a neighborhood of a regular point by the formula
\[
f(x,y) = \begin{cases}
\frac{c_1 x^2}{2c_2}-\frac{1}{2}xy+F\left(c_1x-c_2y\right) & \mbox{if } c_2\neq 0\\
 \frac{1}{2}xy+F(x) & \mbox{if } c_2=0.
 \end{cases}
\]
\end{thm}

\begin{proof}[Proof of Theorem \ref{main}]
It follows from Corollary \ref{constantgcor} that
\[
c_1(x-2f_y(x,y))=c_2(y+2f_x(x,y))
\]
where $(c_1,c_2)\neq(0,0)$ is a pair of constants.

By using the method of characteristics (see \cite{Ev}), we obtain
\[
f(x,y) = \begin{cases}
\frac{c_1 x^2}{2c_2}-\frac{1}{2}xy+F\left(c_1x-c_2y\right) & \mbox{if } c_2\neq 0\\
 \frac{1}{2}xy+F(x) & \mbox{if } c_2=0
 \end{cases}
\]
for some function $F$.

\end{proof}

\smallskip

\section{On the p-minimal graph equation}

In this section, we show that the p-minimal graph equation also splits into two first order PDEs. Moreover, a formula of the solution near a regular point can be written down explicitly.

Recall that the p-minimal graph equation is given by
\begin{equation}\label{pminimal}
(x-2f_y)^2f_{xx}+2(x-2f_y)(y+2f_x)f_{xy}+(y+2f_x)^2f_{yy}=0.
\end{equation}
Note that the equation above coincides with that in \cite{ChHwMaYa1} if we set $u=-2f$.

A computation gives the following.

\begin{prop}
Let $U$ be an open set where $x-2f_y\neq 0$ (resp. $y+2f_x\neq 0$). Let $f$ be a $C^2$ solution of the equation (\ref{pminimal}) on $U$. Then
\[
g(x,y)=\frac{y+2f_x}{x-2f_y}\quad \left(\text{resp. } h(x,y)=\frac{x-2f_y}{y+2f_x}\right)
\]
is a solution of the inviscid Burgers' equation
\[
g_x=-gg_y \quad \left(\text{resp. } h_y=-hh_x\right).
\]
\end{prop}

Next, we give a formula to the local solution of (\ref{pminimal}) near a regular point $(x_0,y_0)$. Only the case $x_0-2f_y(x_0,y_0)\neq 0$ will be considered. The case $y_0+2f_x(x_0,y_0)\neq 0$, being very similar, will be omitted. 

\begin{thm}
Let $f$ be a $C^2$ solution of (\ref{pminimal}). Assume that $x_0-2f_y(x_0,y_0)\neq 0$. Then the following holds in a neighborhood of $(x_0,y_0)$
\[
\begin{split}
&f(x,y)=\frac{1}{2}(-\tilde y(x,y)+x_0F(\tilde y(x,y)))(x-x_0)+G(\tilde y(x,y))\\
&y=(x-x_0)F(\tilde y(x,y))+\tilde y(x,y)
\end{split}
\]
for some $C^2$ function $F,G:\Real\to\Real$.
\end{thm}

Before giving the proof of the theorem, let us recover the global solution (\ref{Bern}) from the above formulas for $f$. For simplicity, we assume that $x_0=0$. If we let $F$ be the constant function $F\equiv c$, then 
\[
f(x,y)=\frac{1}{2}(-yx+cx^2)+G(y-cx)
\]
which is the same as the second family in (\ref{Bern}) (with $a\neq 0$). 

If we set $F$ and $G$ to linear functions $F(r)=c_1r+c_0$ and $G(r)=d_1r+d_0$, then 
\[
f(x,y)=\left(-\frac{1}{2}x+d_1\right)\frac{y-c_0x}{c_1x+1}+d_0. 
\]
If we set $d_1=-\frac{1}{2c_1}$, then we obtain 
\[
f(x,y)=d_1y-dc_0x+d_0
\]
which is the same as the first family in (\ref{Bern}) (with $b\neq 0$). 

\begin{proof}
Let us fix a regular point $(x_0,y_0)$ and consider the Hamiltonian system of the Hamiltonian $H_1(x,y,z,p,q)=p+zq$
\[
\dot x=1, \quad \dot y=z, \quad \dot p=-pq, \quad  \dot q=-q^2,\quad \dot z(t)=p+zq
\]
with initial conditions $x(0)=x_0$, $y(0)=\bar y$, $p(0)=-F(\bar y)F'(\bar y)$, $q(0)=F'(\bar y)$, and $z(0)=F(\bar y)$.

The solution is given by
\[
\begin{split}
&x(t)=t+x_0,\quad y(t)=tF(\bar y)+\bar y,\quad p(t)=-\frac{F(\bar y)F'(\bar y)}{1+tF'(\bar y)},\\
&q(t)=\frac{F'(\bar y)}{1+tF'(\bar y)},\quad z(t)=F(\bar y).
\end{split}
\]

Therefore, by the method of characteristics, the solution $g$ is given by
\[
g(x,y)=F(\bar y(x,y)),\quad y=(x-x_0)F(\bar y(x,y))+\bar y(x,y).
\]

Let $H_2(x,y,z,p,q)=\frac{1}{2}\left(y+2p\right)-\frac{1}{2}g(x,y)(x-2q)$ and consider the corresponding Hamiltonian system
\[
\dot x=1,\quad \dot y=g(x,y),\quad \dot p=\frac{1}{2}\left(g_x(x,y)(x-2q)+g(x,y)\right),
\]
\[
\dot q=\frac{1}{2}\left(g_y(x,y)(x-2q)-1\right),\quad \dot z=p+qg(x,y)
\]
with initial conditions $x(0)=x_0$, $y(0)=\tilde y$, $q(0)=G'(\tilde y)$, $p(0)=\frac{1}{2}\left(-\tilde y+g(x_0,\tilde y)(x_0-2G'(\tilde y)\right))$, and $z(0)=G(\tilde y)$.

A computation shows that $g(x(t),y(t))$ is independent of $t$. Therefore,
\[
x(t)=t+x_0,\quad y(t)=tg(x_0,\tilde y)+\tilde y=tF(\tilde y)+\tilde y.
\]

Moreover, we have
\[
\ddot z=\frac{d}{dt}\left(p+\frac{1}{2}qg(x,y)\right)=0.
\]

Since
\[
\dot z(0)=p(0)+q(0)g(x(0),y(0))=\frac{1}{2}\left(-\tilde y+x_0F(\tilde y)\right),
\]
we also have
\[
z(t)=\frac{1}{2}\left(-\tilde y+x_0F(\tilde y)\right)t+G(\tilde y).
\]

Therefore, the solution $f$ is given by
\[
\begin{split}
&f(x,y)=\frac{1}{2}(-\tilde y(x,y)+x_0F(\tilde y(x,y)))(x-x_0)+G(\tilde y(x,y))\\
&y=(x-x_0)F(\tilde y(x,y))+\tilde y(x,y).
\end{split}
\]

\end{proof}

\end{document}